\documentclass[a4paper,twoside]{article}
\usepackage{a4}
\usepackage{amssymb}
\usepackage{amsmath}
\usepackage{upref}
\usepackage[active]{srcltx}
\usepackage[pagebackref,colorlinks,citecolor=blue,linkcolor=blue]{hyperref}
\allowdisplaybreaks[2] 
\newcount\minutes \newcount\hours
\hours=\time \divide\hours 60 \minutes=\hours
 \multiply\minutes-60
\advance\minutes \time
\newcommand{\klockan}{\the\hours:{\ifnum\minutes<10 0\fi}\the\minutes}
\newcommand{\tid}{\today\ \klockan}
\newcommand{\prtid}{\smash{\raise 10mm \hbox{\LaTeX ed \tid}}}
\renewcommand{\prtid}{}
\makeatletter  \headheight 10pt
\def\sectionmark#1{} 
\def\subsectionmark#1{}
\newcommand{\sectnr}{\ifnum \c@secnumdepth >\z@
                 \thesection.\hskip 1em\relax \fi}
\def\@evenhead{\footnotesize\rm\thepage\hfil\leftmark\hfil\llap{\prtid}}
\def\@oddhead{\footnotesize\rm\rlap{\prtid}\hfil\rightmark\hfil\thepage}
\def\tableofcontents{\section*{Contents} 
 \@starttoc{toc}}
\makeatother
\makeatletter
\def\@biblabel#1{#1.}
\makeatother
\makeatletter
\let\Thebibliography=\thebibliography
\renewcommand{\thebibliography}[1]{\def\@mkboth##1##2{}\Thebibliography{#1}
\addcontentsline{toc}{section}{References}
\frenchspacing 
\setlength{\@topsep}{0pt}
\setlength{\itemsep}{0pt}%
\setlength{\parskip}{0pt plus 2pt}%
} \makeatother
\makeatletter
\def\mdots@{\mathinner.\nonscript\!.%
 \ifx\next,.\else\ifx\next;.\else\ifx\next..\else
 \nonscript\!\mathinner.\fi\fi\fi}
\let\ldots\mdots@
\let\cdots\mdots@
\let\dotso\mdots@
\let\dotsb\mdots@
\let\dotsm\mdots@
\let\dotsc\mdots@
\def\vdots{\vbox{\baselineskip2.8\p@ \lineskiplimit\z@
    \kern6\p@\hbox{.}\hbox{.}\hbox{.}\kern3\p@}}
\def\ddots{\mathinner{\mkern1mu\raise8.6\p@\vbox{\kern7\p@\hbox{.}}%
    \raise5.8\p@\hbox{.}\raise3\p@\hbox{.}\mkern1mu}}
\makeatother
\makeatletter
\let\Enumerate=\enumerate
\renewcommand{\enumerate}{\Enumerate%
\setlength{\@topsep}{0pt}
\setlength{\itemsep}{0pt}%
\setlength{\parskip}{0pt plus 1pt}%
\renewcommand{\theenumi}{\textup{(\alph{enumi})}}%
\renewcommand{\labelenumi}{\theenumi}%
}
\let\endEnumerate=\endenumerate
\renewcommand{\endenumerate}{\endEnumerate\unskip}
\makeatother
\makeatletter
\def\@seccntformat#1{\csname the#1\endcsname.\quad}
\makeatother
\usepackage{ifthen}
\newcommand{\authortitle}[3]{\author{#1}\markboth{#1}{#2}\ifthenelse{\equal{#3}{}}{\title{#2}}{\title{#3}}}
\newcommand{\art}[6]{{\sc #1, \rm #2, \it #3 \bf #4 \rm (#5), \mbox{#6}.}}
\newcommand{\auth}[2]{{#1, #2.}}

\newcommand{\book}[3]{{\sc #1, \it #2, \rm #3.}}
\newcommand{\AND}{{\rm and }}
\RequirePackage{amsthm}
\newtheoremstyle{descriptive}%
  {\topsep}   
  {\topsep}   
  {\rmfamily} 
  {}          
  {\bfseries} 
  {.}         
  { }         
  {}          
\newtheoremstyle{propositional}%
  {\topsep}   
  {\topsep}   
  {\itshape}  
  {}          
  {\bfseries} 
  {.}         
  { }         
  {}          
\theoremstyle{propositional}
\newtheorem{thm}{Theorem}[section]
\newtheorem{theorem}[thm]{Theorem}  
\newtheorem{proposition}[thm]{Proposition}
\theoremstyle{descriptive}
\newtheorem{definition}[thm]{Definition}
\newtheorem{example}[thm]{Example}
\makeatletter
\renewenvironment{proof}[1][\proofname]{\par
  \pushQED{\qed}%
  \normalfont
  \trivlist
  \item[\hskip\labelsep
        \itshape
    #1\@addpunct{.}]\ignorespaces
}{%
  \popQED\endtrivlist\@endpefalse
} \makeatother
\newcommand{\setm}{\setminus}
\renewcommand{\emptyset}{\varnothing}
\def\vint{\mathop{\mathchoice%
          {\setbox0\hbox{$\displaystyle\intop$}\kern 0.22\wd0%
           \vcenter{\hrule width 0.6\wd0}\kern -0.82\wd0}%
          {\setbox0\hbox{$\textstyle\intop$}\kern 0.2\wd0%
           \vcenter{\hrule width 0.6\wd0}\kern -0.8\wd0}%
          {\setbox0\hbox{$\scriptstyle\intop$}\kern 0.2\wd0%
           \vcenter{\hrule width 0.6\wd0}\kern -0.8\wd0}%
          {\setbox0\hbox{$\scriptscriptstyle\intop$}\kern 0.2\wd0%
           \vcenter{\hrule width 0.6\wd0}\kern -0.8\wd0}}%
          \mathopen{}\int}
\newcommand{\Cp}{{C_p}}
\DeclareMathOperator{\Div}{div}
\DeclareMathOperator{\diam}{diam}
\DeclareMathOperator{\dist}{dist}
\DeclareMathOperator{\Lip}{Lip}
\newcommand{\bdy}{\partial}
\newcommand{\loc}{_{\rm loc}}
{\catcode`p =12 \catcode`t =12 \gdef\eeaa#1pt{#1}}      
\def\accentadjtext#1{\setbox0\hbox{$#1$}\kern   
                \expandafter\eeaa\the\fontdimen1\textfont1 \ht0 }
\def\accentadjscript#1{\setbox0\hbox{$#1$}\kern 
                \expandafter\eeaa\the\fontdimen1\scriptfont1 \ht0 }
\def\accentadjscriptscript#1{\setbox0\hbox{$#1$}\kern   
                \expandafter\eeaa\the\fontdimen1\scriptscriptfont1 \ht0 }
\def\accentadjtextback#1{\setbox0\hbox{$#1$}\kern       
                -\expandafter\eeaa\the\fontdimen1\textfont1 \ht0 }
\def\accentadjscriptback#1{\setbox0\hbox{$#1$}\kern     
                -\expandafter\eeaa\the\fontdimen1\scriptfont1 \ht0 }
\def\accentadjscriptscriptback#1{\setbox0\hbox{$#1$}\kern 
                -\expandafter\eeaa\the\fontdimen1\scriptscriptfont1 \ht0 }
\def\itoverline#1{{\mathsurround0pt\mathchoice
        {\rlap{$\accentadjtext{\displaystyle #1}
                \accentadjtext{\vrule height1.593pt}
                \overline{\phantom{\displaystyle #1}
                \accentadjtextback{\displaystyle #1}}$}{#1}}
        {\rlap{$\accentadjtext{\textstyle #1}
                \accentadjtext{\vrule height1.593pt}
                \overline{\phantom{\textstyle #1}
                \accentadjtextback{\textstyle #1}}$}{#1}}
        {\rlap{$\accentadjscript{\scriptstyle #1}
                \accentadjscript{\vrule height1.593pt}
                \overline{\phantom{\scriptstyle #1}
                \accentadjscriptback{\scriptstyle #1}}$}{#1}}
        {\rlap{$\accentadjscriptscript{\scriptscriptstyle #1}
                \accentadjscriptscript{\vrule height1.593pt}
                \overline{\phantom{\scriptscriptstyle #1}
                \accentadjscriptscriptback{\scriptscriptstyle #1}}$}{#1}}}}
\def\itunderline#1{{\mathsurround0pt\mathchoice
        {\rlap{$\underline{\phantom{\displaystyle #1}
                \accentadjtextback{\displaystyle #1}}$}{#1}}
        {\rlap{$\underline{\phantom{\textstyle #1}
                \accentadjtextback{\textstyle #1}}$}{#1}}
        {\rlap{$\underline{\phantom{\scriptstyle #1}
                \accentadjscriptback{\scriptstyle #1}}$}{#1}}
        {\rlap{$\underline{\phantom{\scriptscriptstyle #1}
                \accentadjscriptscriptback{\scriptscriptstyle #1}}$}{#1}}}}
\newcommand{\dmu}{d\mu}
\newcommand{\ds}{ds}
\newcommand{\eps}{\varepsilon}
\newcommand{\bdyOmegaX}{\bdy\Omega\setm\{\infty\}}
\newcommand{\clOmega}{\overline\Omega}
\newcommand{\clV}{\overline{V}}

\renewcommand{\phi}{\varphi}
\newcommand{\p}{{$p\mspace{1mu}$}}
\newcommand{\R}{\mathbb{R}}
\newcommand{\eR}{{\overline{\R\kern-0.08em}\kern 0.08em}} 

\newcommand{\limplus}{{\mathchoice{\vcenter{\hbox{$\scriptstyle +$}}}
  {\vcenter{\hbox{$\scriptstyle +$}}}
  {\vcenter{\hbox{$\scriptscriptstyle +$}}}
  {\vcenter{\hbox{$\scriptscriptstyle +$}}}
}}
\newcommand{\Lp}{L^p}
\newcommand{\Lploc}{L^{p}\loc}
\newcommand{\Np}{N^{1,p}}
\newcommand{\Nploc}{N^{1,p}\loc}
\newcommand{\Dp}{D^p}
\newcommand{\Dploc}{D^{p}\loc}
\newcommand{\ft}{\tilde{f}}
\newcommand{\ut}{\tilde{u}}
\newcommand{\uHp}{\itoverline{P}}   
\newcommand{\lHp}{\itunderline{P}}  
\newcommand{\oHp}{H}                
\newcommand{\Hp}{P}                 
\newcommand{\upharm}{{\overline{\omega}}}  
\newcommand{\lpharm}{{\itunderline{\omega}}} 

\RequirePackage{mathrsfs}
\newcommand{\K}{{\mathscr{K}}}
\newcommand{\F}{\mathscr{F}}%
\newcommand{\U}{\mathscr{U}}%
\DeclareMathOperator*{\essliminf}{ess\,lim\,inf}
\DeclareMathOperator*{\essinf}{ess\,inf}
\newcommand{\px}{{\ensuremath{p(\cdot)}}}
\numberwithin{equation}{section}
\newcommand{\eqv}{\ensuremath{
\mathchoice{\quad \Longleftrightarrow \quad}{\Leftrightarrow}
                {\Leftrightarrow}{\Leftrightarrow}} }
\newcommand{\imp}{\ensuremath{\Rightarrow} }

\newenvironment{ack}{\medskip{\it Acknowledgement.}}{}

\newcounter{saveenumi}
\begin{document}

\authortitle{Anders Bj\"orn and Daniel Hansevi}
{Semiregular and strongly irregular boundary points on unbounded sets}
{Semiregular and strongly irregular boundary points for \p-harmonic functions 
on unbounded sets \\ in metric spaces}

\author
{Anders Bj\"orn \\
\it\small Department of Mathematics, Link\"oping University, \\
\it\small SE-581 83 Link\"oping, Sweden\/{\rm ;}
\it \small anders.bjorn@liu.se
\\
\\
Daniel Hansevi \\
\it\small Department of Mathematics, Link\"oping University, \\
\it\small SE-581 83 Link\"oping, Sweden\/{\rm ;}
\it \small daniel.hansevi@liu.se
}

\date{Preliminary version, \today}
\date{}

\maketitle

\noindent{\small
{\bf Abstract}. 
The trichotomy between regular, semiregular, 
and strongly irregular boundary points for \p-harmonic functions 
is obtained for 
unbounded open sets in 
complete metric spaces with a doubling measure 
supporting a \p-Poincar\'e inequality, $1<p<\infty$. 
We show that these are local properties. 
We also deduce 
several characterizations of semiregular points 
and strongly irregular points. In particular,
semiregular points 
are characterized by means of 
capacity, \p-harmonic measures, removability, and semibarriers. 

\bigskip

\noindent {\small \emph{Key words and phrases}:
barrier, 
boundary regularity, 
Dirichlet problem, 
doubling measure, 
metric space, 
nonlinear potential theory, 
Perron solution, 
\p-harmonic function, 
Poincar\'e inequality, 
semibarrier, 
semiregular boundary point,
strongly irregular boundary point.
}

\medskip

}

\section{Introduction}
\label{sec:intro} 
Let $\Omega\subset\R^n$ be a nonempty bounded open set 
and let $f\in C(\bdy\Omega)$. 
The Perron method
provides us with a unique function $\Hp f$ 
that is harmonic in $\Omega$ and takes the boundary values $f$ 
in a weak sense, i.e., 
$\Hp f$ is a solution of the Dirichlet problem for the Laplace equation
$\Delta u=0$. 
It was introduced on $\R^2$ in 1923 by Perron~\cite{Perron23} 
and independently by Remak~\cite{remak}. 
A point $x_0\in\bdy\Omega$ is 
\emph{regular} if 
$\lim_{\Omega\ni y\to x_0}\Hp f(y)=f(x_0)$ for every $f\in C(\bdy\Omega)$. 
Wiener~\cite{Wiener} characterized 
regular boundary points by means of the 
\emph{Wiener criterion} 
in 1924. 
In the same year  Lebesgue~\cite{Lebesgue} 
gave a different characterization using barriers.

This definition of boundary regularity can be paraphrased in the following way: 
The point $x_0\in\bdy\Omega$ is regular if the following two conditions hold:
\begin{enumerate}
\renewcommand{\theenumi}{\textup{(\roman{enumi})}}%
\item For all $f\in C(\bdy\Omega)$ the limit 
$\lim_{\Omega\ni y\to x_0}\Hp f(y)$ exists. 
\item For all $f\in C(\bdy\Omega)$ there is a sequence 
$\Omega\ni y_j\to x_0$ such that $\lim_{j\to\infty}\Hp f(y_j)=f(x_0)$.
\end{enumerate}
Perhaps surprisingly, 
it is the case that for irregular boundary points 
\emph{exactly one} of these two properties fails; 
one might have guessed that 
both can fail at the same time 
but this can in fact never happen. 
A boundary point $x_0\in\bdy\Omega$ is 
\emph{semiregular} 
if the first condition holds but not the second; 
and \emph{strongly irregular} 
if the second condition holds but not the first. 

For the Laplace equation it is well known that all boundary
points are either regular, semiregular, or strongly irregular,
and this trichotomy (in an abstract linear setting) was developed 
in detail in Luke\v{s}--Mal\'y~\cite{lukesmaly}. 
Key examples of semiregular and strongly irregular points 
are Zaremba's punctured ball and 
the Lebesgue spine, respectively,
see Examples~13.3 and~13.4 in \cite{BBbook}.

A nonlinear analogue is to consider the 
Dirichlet problem for \p-harmonic functions, 
which are solutions of the \p-Laplace equation 
$\Delta_p u:=\Div(|\nabla u|^{p-2}\,\nabla u)=0$, $1<p<\infty$. 
This leads to a nonlinear potential theory 
that has been studied since the 1960s. 
Initially, it was developed for $\R^n$, 
but it has also been extended to weighted $\R^n$, 
Riemannian manifolds, 
and other settings. 
In more recent years, it has been generalized to metric spaces, 
see, e.g., the monograph 
Bj\"orn--Bj\"orn~\cite{BBbook} and the references therein. 
The Perron method was extended to such metric spaces by 
Bj\"orn--Bj\"orn--Shanmugalingam~\cite{BBS2} 
for bounded open sets and Hansevi~\cite{hansevi2} for 
unbounded open sets. 

Boundary regularity for \p-harmonic functions on metric spaces 
was first studied by Bj\"orn~\cite{BjIll} and 
Bj\"orn--MacManus--Shan\-mu\-ga\-lin\-gam~\cite{BMS}, 
and a rather extensive study was undertaken 
by Bj\"orn--Bj\"orn~\cite{BB} on bounded open sets. 
Recently this theory was generalized to 
unbounded open sets by Bj\"orn--Hansevi~\cite{BHan1}; 
see also Bj\"orn--Bj\"orn--Li~\cite{BBLi}. 
For further references and 
a historical discussion on regularity for \p-harmonic functions 
we refer the interested reader to the introduction 
in \cite{BHan1}. 

For \p-harmonic functions on $\R^n$ and metric spaces 
the trichotomy was obtained by Bj\"orn~\cite{ABclass} for bounded open sets. 
It was also obtained for unbounded sets in certain Ahlfors regular metric spaces 
by Bj\"orn--Bj\"orn--Li~\cite{BBLi}. 
Adamowicz--Bj\"orn--Bj\"orn~\cite{ABB} obtained 
the trichotomy for $p(\cdot)$-harmonic functions on 
bounded open sets in $\R^n$. 

In this paper we obtain the trichotomy in the following form, 
where regularity is defined using upper Perron solutions 
(Definition~\ref{def:reg}). 
(We use upper Perron solutions as 
it is not known whether  
continuous functions are resolutive 
with respect to unbounded \p-hyperbolic sets.)
\begin{theorem}\label{thm:trichotomy}
\textup{(Trichotomy)} 
Assume that $X$ is a complete metric space equipped with a doubling measure 
supporting a \p-Poincar\'e inequality\textup{,} $1<p<\infty$. 
Let\/ $\Omega\subset X$ be a nonempty\/ 
\textup{(}possibly unbounded\/\textup{)} open set with 
the capacity $\Cp(X\setm\Omega)>0$. 

Let $x_0\in\bdyOmegaX$. 
Then $x_0$ is either 
regular\textup{,} semiregular\textup{,} or strongly irregular 
for functions that are \p-harmonic in $\Omega$. 
Moreover, 
\begin{itemize}
\item $x_0$ is strongly irregular 
if and only if 
$x_0\in\itoverline{R}\setm R$, 
where 
\[
	R
	:= \{x\in\bdyOmegaX:x\text{ is regular}\}.
\]
\item The relatively open set 
\begin{equation} \label{eq-S}
	S
	:= \{x\in\bdyOmegaX:
		\text{there is $r>0$ such that $\Cp(B(x,r)\cap\bdy\Omega)=0$}\}
\end{equation}
consists exactly of all semiregular boundary points of $\bdyOmegaX$.
\end{itemize}
\end{theorem}
The importance of the distinction between semiregular and 
strongly irregular boundary points is perhaps best illustrated by the 
equivalent characterizations 
given in Theorems~\ref{thm:rem-irr-char} 
and~\ref{thm:ess-irr-char}. 
Semiregular points are in some ways not seen by Perron solutions.

Our contribution here is to extend the results in \cite{ABclass} 
to unbounded open sets. 
In order to do so there are extra complications, 
most notably the fact that it is not known whether 
continuous functions are resolutive 
with respect to unbounded \p-hyperbolic sets. 
We will also rely on the recent results by Bj\"orn--Hansevi~\cite{BHan1} 
on regularity for \p-harmonic functions on unbounded sets 
in metric spaces. 
Most of our results are new also on unweighted $\R^n$.

\begin{ack}
The first author was supported by the Swedish Research Council,
grant 2016-03424.
\end{ack}

\section{Notation and preliminaries}
\label{sec:prel} 
We assume that $(X,d,\mu)$ 
is a metric measure space (which we simply refer to as $X$) 
equipped with a metric $d$ and a 
positive complete Borel measure $\mu$ such that 
$0<\mu(B)<\infty$ 
for every ball $B\subset X$. 
It follows that $X$ is second countable. 
For balls 
$B(x_0,r):=\{x\in X:d(x,x_0)<r\}$ and $\lambda>0$, 
we let $\lambda B=\lambda B(x_0,r):=B(x_0,\lambda r)$. 
The $\sigma$-algebra 
on which $\mu$ is defined 
is the completion of the Borel $\sigma$-algebra. 
We also assume that $1<p<\infty$. 
Later we will impose further requirements 
on the space and on the measure. 
We will keep the discussion short, 
see the monographs
Bj\"orn--Bj\"orn~\cite{BBbook} and 
Heinonen--Koskela--Shanmugalingam--Tyson~\cite{HKSTbook} 
for proofs, 
further discussion, 
and references on the topics in this section. 

The measure $\mu$ is 
\emph{doubling} if there exists 
a constant $C\geq 1$ such that 
\[
	0 
	< \mu(2B) 
	\leq C\mu(B) 
	< \infty 
\]
for every ball $B\subset X$. 
A metric space is 
\emph{proper} 
if all bounded closed subsets are compact, 
and this is in particular true if the metric space 
is complete and the measure is doubling. 

We say that a property holds for \emph{\p-almost every curve} 
if it fails only for 
a curve family $\Gamma$ with zero \p-modulus, 
i.e., 
there exists a nonnegative $\rho\in\Lp(X)$ such that 
$\int_\gamma\rho\,\ds=\infty$ for every curve $\gamma\in\Gamma$. 
For us, a curve in $X$ is a rectifiable nonconstant continuous mapping 
from a compact interval into $X$, 
and it can thus be parametrized 
by its arc length $\ds$.

Following Koskela--MacManus~\cite{KoMac98} 
we make the following definition, 
see also Heinonen--Koskela~\cite{HeKo98}.
\begin{definition}\label{def:upper-gradients}
A measurable function $g\colon X\to[0,\infty]$ is a 
\emph{\p-weak upper gradient} 
of the function $u\colon X\to\eR:=[-\infty,\infty]$ 
if 
\[
	|u(\gamma(0)) - u(\gamma(l_{\gamma}))| 
	\leq \int_{\gamma}g\,\ds
\]
for \p-almost every curve 
$\gamma\colon[0,l_{\gamma}]\to X$, 
where we use the convention that the left-hand side is $\infty$ 
whenever at least one of the terms on the left-hand side is infinite. 
\end{definition}
One way of controlling functions by their 
\p-weak upper gradients is to require 
a Poincar\'e inequality to hold.
\begin{definition}\label{def:Poincare-inequality}
We say that $X$ supports a \p-\emph{Poincar\'e inequality} 
if there exist constants, 
$C>0$ and $\lambda\geq 1$ (the dilation constant), 
such that for all balls $B\subset X$, 
all integrable functions $u$ on $X$, 
and all \p-weak upper gradients $g$ of $u$,
\begin{equation}\label{def:Poincare-inequality-ineq}
	\vint_B|u-u_B|\,\dmu
	\leq C\diam(B)\biggl(\vint_{\lambda B}g^p\,\dmu\biggr)^{1/p},
\end{equation}
where $u_B:=\vint_B u\,\dmu:=\frac{1}{\mu(B)}\int_B u\,\dmu$.
\end{definition}
Shanmugalingam~\cite{Shanmugalingam00} 
used \p-weak upper gradients to 
define so-called Newtonian spaces. 
\begin{definition}\label{def:Newtonian-space}
The \emph{Newtonian space} on $X$, 
denoted 
$\Np(X)$, 
is the space of all 
extended real-valued functions $u\in\Lp(X)$ 
such that  
\[
	\|u\|_{\Np(X)} 
	:= \biggl(\int_X|u|^p\,\dmu + \inf_g\int_X g^p\,\dmu\biggr)^{1/p}<\infty, 
\]
where the infimum is taken over all \p-weak upper gradients $g$ of $u$. 
\end{definition}
The quotient space $\Np(X)/\sim$, 
where $u\sim v$ if and only if $\|u-v\|_{\Np(X)}=0$, 
is a Banach space, see 
Shanmugalingam~\cite{Shanmugalingam00}.
\begin{definition}\label{def:Dirichlet-space}
The \emph{Dirichlet space} on $X$, 
denoted 
$\Dp(X)$, 
is the space of all 
measurable extended real-valued functions on $X$ 
that have a \p-weak upper gradient
in $\Lp(X)$. 
\end{definition}
In this paper we assume that functions in $\Np(X)$ and $\Dp(X)$ 
are defined everywhere (with values in $\eR$), 
not just up to an equivalence class. 
This is important, in particular for the definition of 
\p-weak upper gradients to make sense.

A measurable set $A\subset X$ can itself be 
considered to be a metric space 
(with the restriction of $d$ and $\mu$ to $A$) with 
the Newtonian space $\Np(A)$ and the Dirichlet space $\Dp(A)$ 
given by 
Definitions~\ref{def:Newtonian-space}~and~\ref{def:Dirichlet-space}, 
respectively. 
If $X$ is proper and $\Omega\subset X$ is open, 
then $u\in\Nploc(\Omega)$ 
if and only if 
$u\in\Np(V)$ 
for every open $V$ such that $\clV$ 
is a compact subset of $\Omega$, 
and similarly for $\Dploc(\Omega)$. 
If $u\in\Dploc(X)$, 
then there exists 
a \emph{minimal \p-weak upper gradient} $g_u\in\Lploc(X)$ of $u$ such 
that $g_u\leq g$ a.e.\ for 
all \p-weak upper gradients $g\in\Lploc(X)$ of $u$.
\begin{definition}\label{def:capacity}
The (\emph{Sobolev}) \emph{capacity} of a set $E\subset X$ is the number 
\[
	\Cp(E) 
	:= \inf_u\|u\|_{\Np(X)}^p,
\]
where the infimum is taken over all 
$u\in\Np(X)$ such that $u\geq 1$ on $E$. 

A property that holds for all points 
except for those in a set of capacity zero 
is said to hold \emph{quasieverywhere} (\emph{q.e.}). 
\end{definition}
The capacity is countably subadditive, 
and it is the correct gauge 
for distinguishing between two Newtonian functions: 
If $u\in\Np(X)$, then $u\sim v$ if and only if $u=v$ q.e. 
Moreover, 
if $u,v\in\Nploc(X)$ and $u=v$ a.e., then $u=v$ q.e. 

Continuous functions will be assumed to be real-valued 
unless otherwise stated, 
whereas semicontinuous functions are allowed to take values in $\eR$. 
We use the common notation 
$u_\limplus=\max\{u,0\}$, 
let $\chi_E$ denote the characteristic function of the set $E$, 
and consider all neighbourhoods to be open.

\section{The obstacle problem and \texorpdfstring{\boldmath$p\mspace{1mu}$}{p}-harmonic functions}
\label{sec:p-harmonic}
\emph{We assume from now on 
that\/ $1<p<\infty$\textup{,} 
that $X$ is a complete metric measure space 
supporting a \p-Poincar\'e inequality\textup{,} 
that $\mu$ is doubling\textup{,} 
and that\/ $\Omega\subset X$ 
is a nonempty \textup{(}possibly unbounded\textup{)} 
open subset with $\Cp(X\setm\Omega)>0$.}

\medskip

\begin{definition}\label{def:min}
A function $u\in\Nploc(\Omega)$ 
is a \emph{minimizer} 
in $\Omega$ if 
\[
	\int_{\phi\neq 0}g_u^p\,\dmu 
	\leq \int_{\phi\neq 0}g_{u+\phi}^p\,\dmu
	\quad\text{for all }\phi\in\Np_0(\Omega),
\]
where 
$\Np_0(\Omega)=\{u|_\Omega:u\in\Np(X)\text{ and }u=0\text{ in }X\setm\Omega\}$. 
Moreover, 
a function is \emph{\p-harmonic} 
if it is a continuous minimizer. 
\end{definition}
Kinnunen--Shanmugalingam~\cite[Proposition~3.3 and Theorem~5.2]{KiSh01} 
used De Giorgi's method to 
show that every minimizer 
$u$ has a H\"older continuous representative 
$\ut$ such that $\ut=u$ q.e.
Bj\"orn--Marola~\cite[p.\ 362]{BMarola} obtained the same conclusions 
using Moser iterations. 
See alternatively Theorems~8.13 and 8.14 in \cite{BBbook}. 
Note that $\Nploc(\Omega)=\Dploc(\Omega)$, by 
Proposition~4.14 in \cite{BBbook}. 

The following obstacle problem 
is an important tool. 
In this generality, it was 
considered 
by Hansevi~\cite{hansevi1}.
\begin{definition}\label{def:obst}
Let $V\subset X$ be a nonempty open subset with $\Cp(X\setm V)>0$. 
For $\psi\colon V\to\eR$ and $f\in\Dp(V)$, 
let 
\[
	\K_{\psi,f}(V)
	= \{v\in\Dp(V):v-f\in\Dp_0(V)\textup{ and }v\geq\psi\text{ q.e.\ in }V\}, 
\]
where $\Dp_0(V)=\{u|_V:u\in\Dp(X)\text{ and }u=0\text{ in }X\setm V\}$. 
We say that $u\in\K_{\psi,f}(V)$ is a 
\emph{solution of the }$\K_{\psi,f}(V)$-\emph{obstacle problem 
\textup{(}with obstacle $\psi$ and boundary values $f$\,\textup{)}} 
if  
\[
	\int_V g_u^p\,\dmu 
	\leq \int_V g_v^p\,\dmu
	\quad\textup{for all }v\in\K_{\psi,f}(V).
\]
When $V=\Omega$, 
we usually denote $\K_{\psi,f}(\Omega)$ by $\K_{\psi,f}$.
\end{definition}
The $\K_{\psi,f}$-obstacle problem has a unique 
(up to sets of capacity zero) solution 
whenever $\K_{\psi,f}\neq\emptyset$, see 
Hansevi~\cite[Theorem~3.4]{hansevi1}. 
Furthermore, 
there is 
a unique lsc-regularized 
solution of the $\K_{\psi,f}$-obstacle problem, 
by Theorem~4.1 in~\cite{hansevi1}. 
A function $u$ is \emph{lsc-regularized} if $u=u^*$, 
where the \emph{lsc-regularization} $u^*$ of $u$ is defined by 
\[
	u^*(x) 
	= \essliminf_{y\to x}u(y) 
	:= \lim_{r\to 0}\essinf_{B(x,r)}u.
\]

If $\psi\colon\Omega\to[-\infty,\infty)$ 
is continuous as an extended real-valued function, 
and $\K_{\psi,f}\neq\emptyset$, 
then the lsc-regularized solution of the $\K_{\psi,f}$-obstacle problem 
is continuous, 
by Theorem~4.4 in 
\cite{hansevi1}. 
Hence the following generalization 
of Definition~3.3 in 
Bj\"orn--Bj\"orn--Shanmugalingam~\cite{BBS} 
(and Definition~8.31 in \cite{BBbook}) 
to Dirichlet functions and to unbounded sets makes sense. 
It was first used by 
Hansevi~\cite[Definition~4.6]{hansevi1}.
\begin{definition}\label{def:ext}
Let $V\subset X$ be a nonempty open set with $\Cp(X\setm V)>0$. 
The \emph{\p-harmonic extension} 
$\oHp_V f$ of $f\in\Dp(V)$ to $V$ is the continuous solution 
of the $\K_{-\infty,f}(V)$-obstacle problem. 
When $V=\Omega$, we usually write $\oHp f$ instead of $\oHp_\Omega f$.
\end{definition}
\begin{definition}\label{def:superharm}
A function $u\colon\Omega\to(-\infty,\infty]$ 
is \emph{superharmonic} in $\Omega$ if 
\begin{enumerate}
\renewcommand{\theenumi}{\textup{(\roman{enumi})}}%
\item $u$ is lower semicontinuous; 
\item $u$ is not identically $\infty$ in any component of $\Omega$; 
\item for every nonempty open set 
$V$ such that $\clV$ is a compact subset of $\Omega$ 
and all $v\in\Lip(\clV)$, 
we have $\oHp_{V}v\leq u$ in $V$ 
whenever $v\leq u$ on $\bdy V$. 
\end{enumerate}
A function $u\colon\Omega\to[-\infty,\infty)$ is 
\emph{subharmonic} if $-u$ is superharmonic.
\end{definition}
There are several other equivalent definitions of superharmonic functions, 
see, e.g., 
Theorem~6.1 in Bj\"orn~\cite{ABsuper} 
(or Theorem~9.24 and Propositions~9.25 and~9.26 in \cite{BBbook}). 

An lsc-regularized 
solution of the obstacle problem is always superharmonic, 
by Proposition~3.9 in \cite{hansevi1} 
together with Proposition~7.4 in Kinnunen--Martio~\cite{KiMa02} 
(or Proposition~9.4 in \cite{BBbook}). 
On the other hand, 
superharmonic functions are always lsc-regularized, 
by Theorem~7.14 in Kinnunen--Martio~\cite{KiMa02} 
(or Theorem~9.12 in \cite{BBbook}). 

\section{Perron solutions}
\label{sec:perron} 
\emph{In addition to the assumptions given at the 
beginning of Section~\ref{sec:p-harmonic}\textup{,} 
from now on we make the convention that if\/ $\Omega$ is unbounded\textup{,} 
then the point at infinity\textup{,} $\infty$\textup{,} 
belongs to the boundary $\bdy\Omega$. 
Topological notions should therefore be understood with respect to the 
one-point compactification $X^*:=X\cup\{\infty\}$.} 

\medskip

Note that this convention does not affect any of the definitions in Sections~\ref{sec:prel} or~\ref{sec:p-harmonic}, as $\infty$ is \emph{not} added to $X$ (it is added solely to $\bdy\Omega$).

Since continuous functions are assumed to be real-valued, 
every function in $C(\bdy\Omega)$ is bounded 
even if $\Omega$ is unbounded.
Note that since $X$ is second countable so is $X^*$,
and hence $X^*$ is metrizable by Urysohn's metrization theorem,
see, e.g., Munkres~\cite[Theorems~32.3 and~34.1]{Munkres00}.

We will only consider Perron solutions 
and \p-harmonic measures with respect to $\Omega$ 
and therefore omit $\Omega$ from the notation below.
\begin{definition}\label{def:Perron}
Given a function $f\colon\bdy\Omega\to\eR$, 
let $\U_f$ be the collection of all 
functions 
$u$ that are superharmonic in $\Omega$, 
bounded from below, and such that 
\[
	\liminf_{\Omega\ni y\to x}u(y)
	\geq f(x)
	\quad\textup{for all }x\in\bdy\Omega.
\]
The \emph{upper Perron solution} of $f$ is defined by 
\[
	\uHp f(x)
	= \inf_{u\in \U_f }u(x),
	\quad x\in\Omega.
\]
The \emph{lower Perron solution} 
can be defined similarly using subharmonic functions, 
or by letting $\lHp f=-\uHp(-f)$. 
If $\uHp f=\lHp f$, 
then we denote the common value by $\Hp f$. 
Moreover, if $\Hp f$ is real-valued, 
then $f$ is said to be 
\emph{resolutive} (with respect to $\Omega$). 
\end{definition}
An immediate consequence of the definition is that 
$\uHp f\leq\uHp h$ whenever $f\leq h$ on $\bdy\Omega$. 
Moreover, if $\alpha\in\R$ and $\beta\geq 0$, 
then $\uHp(\alpha + \beta f)=\alpha+\beta\uHp f$. 
Corollary~6.3 in Hansevi~\cite{hansevi2} shows that $\lHp f\leq\uHp f$. 
In each component of $\Omega$, $\uHp f$ is either \p-harmonic or 
identically $\pm\infty$, 
by Theorem~4.1 in Bj\"orn--Bj\"orn--Shanmugalingam~\cite{BBS2}
(or Theorem~10.10 in \cite{BBbook}); 
the proof is local and 
applies also to unbounded $\Omega$. 
\begin{definition}\label{def:p-para}
Assume that $\Omega$ is unbounded. 
Then 
$\Omega$ is \emph{\p-parabolic} if 
for every compact $K\subset\Omega$, 
there exist functions $u_j\in\Np(\Omega)$ such that 
$u_j\geq 1$ on $K$ for all $j=1,2,\ldots$\,, and 
\[
	\int_\Omega g_{u_j}^p\,\dmu
	\to 0
	\quad\text{as }j\to\infty.
\]
Otherwise, $\Omega$ is 
\emph{\p-hyperbolic}.
\end{definition}
For examples of \p-parabolic sets, see, e.g., 
Hansevi~\cite{hansevi2}. 
The main reason for introducing \p-parabolic sets 
in \cite{hansevi2} 
was to be able to obtain resolutivity results, 
and in particular, 
establishing the following resolutivity and invariance result 
for \p-parabolic unbounded sets. 
The first such invariance result for \p-harmonic functions 
was obtained, for bounded sets, by Bj\"orn--Bj\"orn--Shanmugalingam~\cite{BBS2}. 
\begin{theorem}\label{thm-hansevi2-main}
\textup{(\cite[Theorem~6.1]{BBS2} 
and~\cite[Theorem~7.8]{hansevi2})}
Assume that\/ $\Omega$ is bounded or \p-parabolic. 
Let $h\colon\bdy\Omega\to\eR$ 
be $0$ q.e.\ on $\bdyOmegaX$ and 
$f\in C(\bdy\Omega)$. 
Then $f$ and $f+h$ are resolutive and $\Hp(f+h)=\Hp f$. 
\end{theorem}
Resolutivity of 
continuous functions is not known for unbounded \p-hyperbolic sets, but it is rather 
trivial to show that constant functions are resolutive. 
We shall 
show that a similar invariance result as in 
Theorem~\ref{thm-hansevi2-main} can be obtained 
for constant functions on unbounded \p-hyperbolic sets. 
This fact will be an important tool when 
characterizing semiregular boundary points. 

We first need to define \p-harmonic measures, 
which despite the name are (usually) 
not 
measures, 
but 
nonlinear generalizations of the harmonic measure.
\begin{definition}\label{def:p-harmonic-measure}
The 
\emph{upper and lower \p-harmonic measures} 
of $E\subset\bdy\Omega$ are 
\[
	\upharm(E)
	:= \uHp\chi_E
	\quad\text{and}\quad
	\lpharm(E) 
	:= \lHp\chi_E,
\]
respectively. 
\end{definition}
\begin{proposition}\label{prop-inv-pharm}
Let $E\subset\bdyOmegaX$\textup{,} 
$a\in\R$\textup{,} and $f\colon\bdy\Omega\to\eR$ be such that 
$\Cp(E)=0$ and $f(x)=a$ for all $x\in\bdy\Omega\setm E$. 
Then $\Hp f\equiv a$. 

In particular\textup{,} $\upharm(E)=\lpharm(E)\equiv 0$. 
\end{proposition}
\begin{proof}
Without loss of generality we may assume that $a=0$. 
As the capacity $\Cp$ is an outer capacity, 
by Corollary~1.3 in Bj\"orn--Bj\"orn--Shanmugalingam~\cite{BBS5} 
(or \cite[Theorem~5.31]{BBbook}), 
we can find open sets $G'_j\supset E$ such that $\Cp(G'_j)<2^{-j-1}$, 
$j=1,2,\ldots$\,. 
From the decreasing sequence 
$\{\bigcup_{k=j}^\infty G'_k\}_{j=1}^\infty$, we can choose
a decreasing subsequence of open sets 
$G_k$ with $\Cp(G_k)<2^{-kp}$, $k=1,2,\ldots$\,. 
By Lemma~5.3 in Bj\"orn--Bj\"orn--Shanmugalingam~\cite{BBS2} 
(or \cite[Lemma~10.17]{BBbook}), 
there is a decreasing sequence 
$\{\psi_j\}_{j=1}^\infty$ of nonnegative functions 
such that $\lim_{j\to\infty}\|\psi_j\|_{\Np(X)}=0$ and 
$\psi_j\geq k-j$ in $G_k$ whenever $k>j$. 
In particular, $\psi_j=\infty$ on $E$ for each $j=1,2,\ldots$\,. 

Let $u_j$ be the lsc-regularized 
solution of the $\K_{\psi_j,0}(\Omega)$-obstacle 
problem, $j=1,2,\ldots$\,. 
As $u_j$ is lsc-regularized 
and $u_j\geq\psi_j$ q.e., 
we see that $u_j\geq k-j$ everywhere in $G_k$ whenever $k>j$, 
and also that $u_j\geq 0$ everywhere in $\Omega$. 
In particular, 
$\liminf_{\Omega\ni y\to x}u_j(y)=\infty$ for $x\in E$, 
which shows that $u_j\in\U_f(\Omega)$ and thus $u_j\geq\uHp f$. 

On the other hand, Theorem~3.2 in Hansevi~\cite{hansevi2} 
shows that the sequence $u_j$ decreases q.e.\ to $0$, 
and hence $\uHp f\leq 0$ q.e.\ in $\Omega$. 
Since $\uHp f$ is continuous, we get that $\uHp f\leq 0$ 
everywhere in $\Omega$. 
Applying this to $-f$ shows that 
$\lHp f=-\uHp(-f)\geq 0$ everywhere in $\Omega$, 
which together with the inequality $\lHp f\leq\uHp f$ shows that 
$\lHp f=\uHp f\equiv 0$. 
In particular, $\lpharm(E)=\lHp\chi_E\equiv 0$ 
and $\upharm(E)=\uHp\chi_E\equiv 0$. 
\end{proof}
We will also need the following result. 
\begin{proposition}\label{prop:Perron-semicont}
If $f\colon\bdy\Omega\to[-\infty,\infty)$ 
is an upper semicontinuous function\textup{,} 
then
\[
	\uHp f
	= \inf_{C(\bdy\Omega)\ni\phi\geq f}\uHp\phi.
\]
\end{proposition}
\begin{proof}
Let $\F=\{\phi\in C(\bdy\Omega):\phi\geq f\}$. 
Then $\F$ is downward directed, 
i.e., for each pair of functions $u,v \in \F$ 
there is a function $w\in\F$ such that $w\leq\min\{u,v\}$. 
Because $f$ is upper semicontinuous, $\bdy\Omega$ is compact, 
and $X^*$ is metrizable, 
it follows from Proposition~1.12 in \cite{BBbook} 
that $f=\inf_{\phi\in\F}\phi$. 
Hence by Lemma~10.31 in~\cite{BBbook} 
(whose proof is valid also for unbounded $\Omega$) 
$\uHp f=\inf_{\phi\in\F}\uHp\phi$.
\end{proof}

\section{Boundary regularity}
\label{sec:bdy-regularity}
It is not known whether 
continuous functions are resolutive also 
with respect to unbounded \p-hyperbolic sets. 
We therefore define regular boundary points in the following way. 
\begin{definition}\label{def:reg}
We say that a boundary point $x_0\in\bdy\Omega$ is \emph{regular} if 
\[
	\lim_{\Omega\ni y\to x_0}\uHp f(y)
	= f(x_0)
	\quad\text{for all }f\in C(\bdy\Omega).
\]
This can be paraphrased in the following way: 
A point $x_0\in\bdy\Omega$ is regular if the following two conditions hold: 
\begin{enumerate}
\renewcommand{\theenumi}{\textup{(\Roman{enumi})}}%
\item\label{semi}
For all $f\in C(\bdy\Omega)$ the limit 
\[
	\lim_{\Omega\ni y\to x_0}\uHp f(y)
	\quad\text{exists}.
\]
\item\label{strong}
For all $f\in C(\bdy\Omega)$ there is a sequence 
$\{y_j\}_{j=1}^\infty$ in $\Omega$ such that 
\[
	\lim_{j\to\infty}y_j
	= x_0
	\quad\text{and}\quad
	\lim_{j\to\infty}\uHp f(y_j)
	= f(x_0).
\]
\end{enumerate}
\bigskip

Furthermore, we say that a boundary point 
$x_0\in\bdy\Omega$ is \emph{semiregular} 
if \ref{semi} holds but not \ref{strong}; 
and \emph{strongly irregular} 
if \ref{strong} holds but not \ref{semi}.
\end{definition}
We do not require $\Omega$ to be bounded in this definition, 
but if it is, 
then it follows from 
Theorem~6.1 in Bj\"orn--Bj\"orn--Shanmugalingam~\cite{BBS2} 
(or Theorem~10.22 in \cite{BBbook}) that our 
definition coincides with the definitions of regularity in 
Bj\"orn--Bj\"orn--Shanmugalingam~\cite{BBS}, \cite{BBS2},
and 
Bj\"orn--Bj\"orn~\cite{BB}, \cite{BBbook}, 
where 
regularity is defined using 
$\Hp f$ or 
$\oHp f$.
Thus we can use the boundary regularity results 
from these papers 
when considering bounded sets.

Since $\uHp f=-\lHp(-f)$, 
the same concept of regularity is obtained if we replace the 
upper Perron solution by the lower Perron solution in 
Definition~\ref{def:reg}. 

Boundary regularity for \p-harmonic functions 
on unbounded sets in metric spaces was 
recently studied by Bj\"orn--Hansevi~\cite{BHan1}. 
We will need some of the characterizations 
obtained therein. 
For the reader's convenience we state these results here. 
We will not discuss regularity of the point $\infty$ 
in this paper. 
One of the important results we will need from \cite{BHan1} 
is the Kellogg property. 
\begin{theorem}\label{thm:kellogg}
\textup{(The Kellogg property)}
If $I$ is the set of irregular points in $\bdyOmegaX$\textup{,} 
then $\Cp(I)=0$.
\end{theorem}
\begin{definition}\label{def:barrier}
A function $u$ is a \emph{barrier} (with respect to $\Omega$) 
at $x_0\in\bdy\Omega$ if 
\begin{enumerate}
\renewcommand{\theenumi}{\textup{(\roman{enumi})}}%
\item\label{barrier-i}
$u$ is superharmonic in $\Omega$;
\item\label{barrier-ii}
$\lim_{\Omega\ni y\to x_0}u(y)=0$;
\item\label{barrier-iii}
$\liminf_{\Omega\ni y\to x}u(y)>0$ for every $x\in\bdy\Omega\setm\{x_0\}$.
\end{enumerate}
\end{definition}
Superharmonic functions satisfy the 
strong minimum principle, i.e., 
if $u$ is superharmonic and attains its minimum in some component 
$G$ of $\Omega$, then $u|_G$ is constant 
(see Theorem~9.13 in \cite{BBbook}). 
This implies that a barrier is always nonnegative, 
and furthermore, that a barrier is positive if 
$\bdy G\setm\{x_0\}\neq\emptyset$ for 
every component $G\subset\Omega$.

The following result is a collection 
of the key facts we will need from 
Bj\"orn--Hansevi~\cite[Theorems~5.2, 5.3, 6.2, and 9.1]{BHan1}. 
\begin{theorem}\label{thm:reg}
Let $x_0\in\bdyOmegaX$ 
and $\delta>0$. 
Also define $d_{x_0}\colon X^*\to[0,1]$ 
by 
\begin{equation}\label{eq-dx0}
	d_{x_0}(x) 
	= \begin{cases}
		\min\{d(x,x_0),1\} 
			& \text{if }x\neq\infty, \\
		1
			& \text{if }x=\infty.
	\end{cases}
\end{equation}
Then the following are equivalent\/\textup{:}
\begin{enumerate}
\item\label{reg-reg} 
The point $x_0$ is regular.
\item\label{barrier-bar-Om}
There is a barrier at $x_0$.
\item\label{barrier-bar-pos-Om}
There is a positive continuous barrier at $x_0$. 
\item\label{barrier-reg-B}
The point $x_0$ is regular with respect to $\Omega\cap B(x_0,\delta)$.
\item\label{reg-cont-x0}
It is true that 
\[
	\lim_{\Omega\ni y\to x_0}\uHp f(y)
	= f(x_0)
\]
for all $f\colon\bdy\Omega\to\R$ that are 
bounded on $\bdy\Omega$ and continuous at $x_0$. 
\item \label{reg-Pd}
It is true that 
\[
	\lim_{\Omega\ni y\to x_0}\uHp d_{x_0}(y)
	= 0.
\]
\item\label{reg-2-obst-dist}
The continuous solution $u$ of the $\K_{d_{x_0},d_{x_0}}$-obstacle 
problem\textup{,} 
satisfies 
\[
	\lim_{\Omega\ni y\to x_0}u(y)
	= 0.
\]
\item\label{reg-2-obst-cont}
If $f\in C(\overline{\Omega})\cap\Dp(\Omega)$\textup{,} 
then the continuous solution $u$ of the $\K_{f,f}$-obstacle 
problem\textup{,} 
satisfies 
\[
	\lim_{\Omega\ni y\to x_0}u(y)
	= f(x_0).
\]
\end{enumerate}
\end{theorem}

\section{Semiregular and strongly irregular points}
\label{sec:trichotomy}
We are now ready to start our discussion of semiregular 
and strongly irregular boundary points. 
We begin by proving Theorem~\ref{thm:trichotomy}.
\begin{proof}[Proof of Theorem~\ref{thm:trichotomy}]
We consider two complementary cases.

\emph{Case} 1: 
\emph{There exists $r>0$ such that $\Cp(B\cap\bdy\Omega)=0$\textup{,} 
where $B:=B(x_0,r)$.} 

Let $G$ be the component of $B$ containing $x_0$. 
Since $X$ is quasiconvex, by, e.g., Theorem~4.32 in \cite{BBbook},
and thus locally connected, it follows that $G$ is open. 
Let $F=G\setm\Omega$. 
Then 
\[
	\Cp(G\cap\bdy F)
	= \Cp(G\cap\bdy\Omega)
	\leq \Cp(B\cap\bdy\Omega)
	= 0,
\]
and hence 
$\Cp(F)=0$, 
by 
Lemma~8.6
in Bj\"orn--Bj\"orn--Shanmugalingam~\cite{BBS2}
(or Lemma~4.5 in \cite{BBbook}). 

Let $f\in C(\bdy\Omega)$. 
Then the Perron solution $\uHp f$ is bounded (as $f$ is bounded), 
and thus $\uHp f$ has a \p-harmonic extension $U$ 
to $\Omega\cup G$, 
by 
Theorem~6.2 in Bj\"orn~\cite{ABremove} 
(or Theorem~12.2 in \cite{BBbook}). 
Since $U$ is continuous, 
it follows that 
\[
	\lim_{\Omega\ni y\to x_0}\uHp f(y)
	= \lim_{\Omega\ni y\to x_0}U(y)
	= U(x_0),
\]
i.e., 
condition~\ref{semi} in Definition~\ref{def:reg}
holds, 
and hence $x_0$ is either regular or semiregular. 

To show that $x_0$ must be semiregular,
we let $f(x)=(1-d_{x_0}(x)/{\min\{r,1\}})_\limplus$ 
on $\bdy\Omega$, 
where $d_{x_0}$ is defined by \eqref{eq-dx0}. 
Then $f=0$ q.e.\ on $\bdy\Omega$, 
and Proposition~\ref{prop-inv-pharm} 
shows that $\Hp f\equiv 0$. 
Since 
\[
	\lim_{\Omega\ni y\to x_0}\uHp f(y)
	= 0
	\neq 1
	= f(x_0),
\]
$x_0$ is not regular, 
and hence must be semiregular.

\smallskip

\emph{Case} 2: 
\emph{For all $r>0$\textup{,} $\Cp(B(x_0,r)\cap\bdy\Omega)>0$.} 

For every $j=1,2,\ldots$\,, 
$\Cp(B(x_0,1/j)\cap\bdy\Omega)>0$, 
and by the Kellogg property (Theorem~\ref{thm:kellogg}) 
there exists a regular boundary point 
$x_j\in B(x_0,1/j)\cap\bdy\Omega$. 
(We do not require the $x_j$ to be distinct.) 

Let $f\in C(\bdy\Omega)$. 
Because $x_j$ is regular, 
there is $y_j\in B(x_j,1/j)\cap\Omega$ 
so that $|\uHp f(y_j)-f(x_j)|<1/j$. 
It follows that $y_j\to x_0$ and $\uHp f(y_j)\to f(x_0)$
as $j\to\infty$, 
i.e., 
condition~\ref{strong} in Definition~\ref{def:reg}
holds, 
and hence $x_0$ must be either regular or strongly irregular.

As there are no strongly irregular points in case~1,
it follows that $x_0\in\bdyOmegaX$
is strongly irregular if and only if $x_0\in\itoverline{R}\setm R$, 
where $R:=\{x\in\bdyOmegaX:x\text{ is regular}\}$.
And since there are no semiregular points in case~2,
the set $S$ in \eqref{eq-S}
consists exactly of all semiregular boundary points of $\bdyOmegaX$.
\end{proof}
In fact, in case~2 
it is possible to improve upon the result above. 
The sequence 
$\{y_j\}_{j=1}^\infty$ 
can be chosen independently of $f$, 
see the characterization \ref{not-reg-one-seq} 
in Theorem~\ref{thm:rem-irr-char}. 

We will characterize 
semiregular points by 
a number of equivalent conditions in Theorem~\ref{thm:rem-irr-char}. 
But first we obtain 
the following characterizations 
of relatively open sets of semiregular points. 
\begin{theorem}\label{thm:irr-char-V}
Let $V\subset\bdyOmegaX$ be relatively open. 
Then the following statements are equivalent\/\textup{:}
\begin{enumerate}
\item\label{V-semireg}
The set $V$ consists entirely of semiregular points.
\item\label{V-R}
The set $V$ does not contain any regular point.
\item\label{V-Cp-V-bdy}
The capacity $\Cp(V)=0$. 
\item\label{V-upharm}
The upper \p-harmonic measure 
$\upharm(V)\equiv 0$.
\item\label{V-lpharm}
The lower \p-harmonic measure 
$\lpharm(V)\equiv 0$.
\item\label{V-alt-def-irr-super}
The set\/ $\Omega\cup V$ is open in $X$\textup{,} 
$\Cp(X\setm(\Omega\cup V))>0$\textup{,} 
$\mu(V)=0$\textup{,} 
and every function that is bounded and superharmonic in $\Omega$ 
has a superharmonic extension to $\Omega\cup V$. 
\item\label{V-alt-def-irr}
\setcounter{saveenumi}{\value{enumi}}
The set\/ $\Omega\cup V$ is open in $X$\textup{,} 
$\Cp(X\setm(\Omega\cup V))>0$\textup{,} 
and every function that is bounded and \p-harmonic in $\Omega$ 
has a \p-harmonic extension to $\Omega\cup V$. 
\end{enumerate}

If moreover 
$\Omega$ is bounded or \p-parabolic\textup{,} 
then also the following statement 
is equivalent to the statements above. 
\begin{enumerate}
\setcounter{enumi}{\value{saveenumi}}
\item\label{V-rem-motiv}
For every $f\in C(\bdy\Omega)$\textup{,} 
the Perron solution $\Hp f$ depends only on $f|_{\bdy\Omega\setm V}$ 
\textup{(}i.e., if $f,h\in C(\bdy\Omega)$ and 
$f=h$ on $\bdy\Omega\setm V$\textup{,} 
then $\Hp f\equiv\Hp h$\textup{)}. 
\end{enumerate}
\end{theorem}
Note that there are examples of sets 
with positive capacity and even positive measure 
which are removable for bounded \p-harmonic functions, 
see Section~9 in Bj\"orn~\cite{ABremove} (or \cite[Section~12.3]{BBbook}). 
For superharmonic functions it is not known whether 
such examples exist. 
This motivates the formulations of 
\ref{V-alt-def-irr-super} 
and~\ref{V-alt-def-irr}. 

The following example shows that the condition 
$\Cp(X\setm(\Omega\cup V))>0$ 
cannot be dropped from \ref{V-alt-def-irr}, 
nor from \ref{alt-def-irr} in 
Theorem~\ref{thm:rem-irr-char} below. 
We do not know whether the conditions 
$\Cp(X\setm(\Omega\cup V))>0$ and $\mu(V)=0$ 
can be dropped from \ref{V-alt-def-irr-super}, 
but they are needed for our proof. 
Similarly they are needed in 
\ref{alt-def-irr-super} in 
Theorem~\ref{thm:rem-irr-char} below. 

The condition $\Cp(X\setm(\Omega\cup V))>0$ 
was unfortunately overlooked in Bj\"orn~\cite{ABclass} 
and in Bj\"orn--Bj\"orn~\cite{BBbook}: 
It should be added to 
conditions (d$'$) and (e$'$) in 
\cite[Theorem~3.1]{ABclass}, 
to (h) and (i) in 
\cite[Theorem~3.3]{ABclass}, 
to (f$'$) and (g$'$) in 
\cite[Theorem~13.5]{BBbook}, 
and to (j) and (l) in 
\cite[Theorem~13.10]{BBbook}. 
\begin{example}
Let $X=[0,1]$ be equipped with the Lebesgue measure, 
and let $1<p<\infty$, $\Omega=(0,1]$ and $V=\{0\}$. 
Then $\Cp(V)>0$. 
In this case the \p-harmonic functions on $\Omega$ are just 
the constant functions, and these trivially 
have \p-harmonic extensions to $X$. 
Thus the condition 
$\Cp(X\setm(\Omega\cup V))>0$ 
cannot be dropped from \ref{V-alt-def-irr}. 

On the other hand, the set $V$ is not removable 
for bounded superharmonic functions on $\Omega$, 
see Example~9.1 in Bj\"orn~\cite{ABremove} 
or Example~12.17 in \cite{BBbook}. 
\end{example}
\begin{proof}[Proof of Theorem~\ref{thm:irr-char-V}]
\ref{V-R} $\imp$ \ref{V-Cp-V-bdy} 
This follows from the Kellogg property (Theorem~\ref{thm:kellogg}).

\smallskip

\ref{V-Cp-V-bdy} $\imp$ \ref{V-upharm}
This follows directly from Proposition~\ref{prop-inv-pharm}.

\smallskip

\ref{V-upharm} $\imp$ \ref{V-lpharm}
This is trivial.

\smallskip

\ref{V-lpharm} $\imp$ \ref{V-R} 
Suppose that $x\in V$ is regular. 
Because $\chi_V$ is continuous at $x$, 
this yields a contradiction, 
as it follows from Theorem~\ref{thm:reg} that 
\[
	0
	= \lim_{\Omega\ni y\to x}\lpharm(V)(y)
	= \lim_{\Omega\ni y\to x}\lHp\chi_V(y)
	= -\lim_{\Omega\ni y\to x}\uHp(-\chi_V)(y)
	= \chi_V(x)
	= 1.
\]
Thus $V$ does not contain any regular point. 

\smallskip

\ref{V-Cp-V-bdy} $\imp$ \ref{V-alt-def-irr-super} 
Suppose that $\Cp(V)=0$. 
Then $\Cp(X\setm(\Omega\cup V))=\Cp(X\setm\Omega)>0$
and $\mu(V)=0$. 
Let $x\in V$ and 
let $G$ be a connected neighbourhood of $x$ such that 
$G\cap\bdy\Omega\subset V$. 
Sets of capacity zero cannot separate space, 
by Lemma~4.6 in Bj\"orn--Bj\"orn~\cite{BBbook}, 
and hence $G\setm\bdy\Omega$ must be connected, i.e., 
$G\subset\clOmega$, 
from which it follows that $\Omega\cup V$ is open in $X$. 
The superharmonic extension is now provided by 
Theorem~6.3 in Bj\"orn~\cite{ABremove} 
(or Theorem~12.3 in \cite{BBbook}).

\smallskip

\ref{V-alt-def-irr-super} $\imp$ \ref{V-alt-def-irr} 
Let $u$ be a bounded \p-harmonic function on $\Omega$. 
Then, by assumption, $u$ has a superharmonic extension $U$ to $\Omega\cup V$. 
Moreover, as $-u$ is also bounded and \p-harmonic, 
there is a superharmonic extension $W$ of $-u$ to $\Omega\cup V$. 
Now, as $-W$ is clearly a subharmonic extension of $u$ to $\Omega\cup V$, 
Proposition~6.5 in Bj\"orn~\cite{ABremove} 
(or Proposition~12.5 in \cite{BBbook}) asserts that 
$U=-W$ is \p-harmonic 
(it is here that we use that $\mu(V)=0$). 

\smallskip

\ref{V-alt-def-irr} $\imp$ \ref{V-semireg} 
Let $x_0\in V$. 
Since $\Omega\cup V$ is open in $X$, we see that $V\cap\bdy(\Omega\cup V)=\emptyset$, 
and hence $x_0\notin\bdy(\Omega\cup V)$. 
Let 
\[
	h(x) 
	= \biggl(1-\frac{d_{x_0}(x)}{\min\{
		\dist(x_0,\bdy(\Omega\cup V)),1\}}\biggr)_\limplus, 
	\quad x\in\bdy\Omega,
\]
where $d_{x_0}$ is defined by \eqref{eq-dx0}.
Then $\uHp h$ is bounded and 
has a \p-harmonic extension $U$ to $\Omega\cup V$, 
and hence the 
Kellogg property (Theorem~\ref{thm:kellogg}) implies that 
\begin{equation}\label{eq-U=0}
	\lim_{\Omega\cup V\ni y\to x}U(y)
	= \lim_{\Omega\ni y\to x}\uHp h(y)
	= h(x)
	= 0
	\quad\text{for q.e. }x\in\bdy(\Omega\cup V)\setm\{\infty\}.
\end{equation}
Let $G$ be the component of $\Omega\cup V$ containing $x_0$. 
Then 
\[
	\Cp(X\setm G) 
	\geq \Cp(X\setm(\Omega\cup V)) 
	>0.
\]
It then follows from 
Lemma~4.3 in Bj\"orn--Bj\"orn~\cite{BB} 
(or Lemma~4.5 in \cite{BBbook})
that $\Cp(\bdy G)>0$. 
In particular, it follows from \eqref{eq-U=0} that $U\not\equiv 1$ in $G$, 
and thus, by the strong maximum principle 
(see Corollary~6.4 in Kinnunen--Shanmugalingam~\cite{KiSh01} 
or \cite[Theorem~8.13]{BBbook}), 
that 
$U(x_0)<1$. 
Therefore 
\[
	\lim_{\Omega\ni y\to x_0}\uHp h(y)
	= U(x_0) 
	< 1 
	= h(x_0),
\]
and hence $x_0$ must be irregular.

However, if $f\in C(\bdy\Omega)$, 
then $\uHp f$ has a \p-harmonic extension $W$ to $\Omega\cup V$. 
Since $W$ is continuous in $\Omega\cup V$, 
it follows that 
\[
	\lim_{\Omega\ni y\to x_0}\uHp f(y)
	= W(x_0),
\]
and hence the limit on the left-hand side always exists. 
Thus $x_0$ is semiregular. 

\smallskip

\ref{V-semireg} $\imp$ \ref{V-R} 
This is trivial.

\bigskip

We now assume that $\Omega$ is bounded or \p-parabolic. 

\smallskip

\ref{V-Cp-V-bdy} $\imp$ \ref{V-rem-motiv}
This implication follows from 
Theorem~\ref{thm-hansevi2-main}.

\smallskip

\ref{V-rem-motiv} $\imp$ \ref{V-lpharm}
As $-\chi_V\colon\bdy\Omega\to\R$ is upper semicontinuous, 
it follows from Proposition~\ref{prop:Perron-semicont}, 
and \ref{V-rem-motiv}, 
that 
\[
	0
	\leq \lpharm(V)
	= -\uHp(-\chi_V)
	= -\inf_{\substack{\phi\in C(\bdy\Omega)\\-\chi_V\leq\phi\leq 0}}\uHp\phi \\
	= 0,
\]
and hence $\lpharm(V)=0$.
\end{proof}
\begin{definition}\label{def:semibarrier}
A function $u$ is a \emph{semibarrier} (with respect to $\Omega$) 
at $x_0\in\bdy\Omega$ if 
\begin{enumerate}
\renewcommand{\theenumi}{\textup{(\roman{enumi})}}%
\item\label{semibarrier-i}
$u$ is superharmonic in $\Omega$;
\item\label{semibarrier-ii}
$\liminf_{\Omega\ni y\to x_0}u(y)=0$;
\item\label{semibarrier-iii}
$\liminf_{\Omega\ni y\to x}u(y)>0$ for every $x\in\bdy\Omega\setm\{x_0\}$.
\end{enumerate}

\smallskip

Moreover, we say that $u$ is a 
\emph{weak semibarrier} 
(with respect to $\Omega$) at $x_0\in\bdy\Omega$ 
if $u$ is a positive superharmonic function such that 
\ref{semibarrier-ii} holds. 
\end{definition}
Now we are ready to characterize 
the semiregular points by means 
of capacity, \p-harmonic measures, removable singularities, 
and semibarriers. 
In particular, we show that semiregularity is a local property. 
\begin{theorem}\label{thm:rem-irr-char}
Let $x_0\in\bdyOmegaX$\textup{,} $\delta>0$\textup{,} and 
$d_{x_0}\colon X^*\to[0,1]$ be defined by \eqref{eq-dx0}. 
Then the following statements are equivalent\/\textup{:} 
\begin{enumerate}
\item\label{semireg}
The point $x_0$ is semiregular. 
\item\label{semireg-local}
The point $x_0$ is semiregular with respect to $G:=\Omega\cap B(x_0,\delta)$. 
\item\label{not-reg-one-seq}
There is no sequence $\{y_j\}_{j=1}^\infty$ in $\Omega$ such that 
$y_j\to x_0$ as $j\to\infty$ and 
\[
	\lim_{j\to\infty}\uHp f(y_j)
	= f(x_0)
	\quad\text{for all }f\in C(\bdy\Omega).
\]
\item\label{not-reg}
The point $x_0$ is neither regular nor strongly irregular. 
\item\label{R}
It is true that 
$x_0\notin\overline{\{x\in\bdy\Omega:x\text{ is regular}\}}$. 
\item\label{Cp-V-bdy}
There is a neighbourhood $V$ of $x_0$ such that $\Cp(V\cap\bdy\Omega)=0$. 
\item\label{Cp-V}
There is a neighbourhood $V$ of $x_0$ 
such that $\Cp(V\setm\Omega)=0$. 
\item\label{upharm}
There is a neighbourhood $V$ of $x_0$ 
such that $\upharm(V\cap\bdy\Omega)\equiv 0$.
\item\label{lpharm}
There is a neighbourhood $V$ of $x_0$ 
such that $\lpharm(V\cap\bdy\Omega)\equiv 0$.
\item\label{alt-def-irr}
There is a neighbourhood $V\subset\clOmega$ of $x_0$\textup{,} 
with $\Cp(X\setm(\Omega\cup V))>0$\textup{,} 
such that every function that is bounded and 
\p-harmonic in $\Omega$ has a \p-harmonic extension to 
$\Omega\cup V$.
\item\label{rem-irr}
There is a neighbourhood $V$ of $x_0$ 
such that every function that is bounded and 
\p-harmonic in $\Omega$ has a \p-harmonic extension to 
$\Omega\cup V$\textup{,} and moreover $x_0$ is irregular.
\item\label{alt-def-irr-super}
There is a neighbourhood $V$ of $x_0$\textup{,} 
with $\Cp(X\setm(\Omega\cup V))>0$ 
and $\mu(V\setm\nobreak\Omega)=0$\textup{,} 
such that every function that is bounded and 
superharmonic in $\Omega$ has a superharmonic extension to 
$\Omega\cup V$.
\item\label{d-lim}
It is true that 
\[
	\lim_{\Omega\ni y\to x_0}\uHp d_{x_0}(y)
	> 0.
\] 
\item\label{d-liminf}
It is true that 
\[
	\liminf_{\Omega\ni y\to x_0}\uHp d_{x_0}(y)
	> 0.
\] 
\item\label{weaksemibarrier}
There is no weak semibarrier at $x_0$. 
\item\label{semibarrier}
There is no semibarrier at $x_0$. 
\item\label{obst-dist-semibarrier}
\setcounter{saveenumi}{\value{enumi}}
The continuous solution of the $\K_{d_{x_0},d_{x_0}}$-obstacle problem 
is not a semibarrier at $x_0$. 
\end{enumerate}

If moreover 
$\Omega$ is bounded or \p-parabolic\textup{,} 
then also the following statement 
is equivalent to the statements above. 
\begin{enumerate}
\setcounter{enumi}{\value{saveenumi}}
\item \label{rem-motiv}
There is a neighbourhood $V$ of $x_0$ such that 
for every $f\in C(\bdy\Omega)$\textup{,} 
the Perron solution $\Hp f$ depends only on $f|_{\bdy\Omega\setm V}$ 
\textup{(}i.e., if $f,h\in C(\bdy\Omega)$ and 
$f=h$ on $\bdy\Omega\setm V$\textup{,} 
then $\Hp f\equiv\Hp h$\textup{)}. 
\end{enumerate}
\end{theorem}
\begin{proof}
\ref{R} 
$\eqv$ \ref{Cp-V-bdy} 
$\eqv$ \ref{upharm} 
$\eqv$ \ref{lpharm} 
$\imp$ \ref{semireg} 
This follows directly from Theorem~\ref{thm:irr-char-V}, 
with $V$ therein 
corresponding to $V\cap\bdy\Omega$ here. 

\smallskip

\ref{semireg} $\imp$ \ref{d-lim} 
Since $x_0$ is semiregular, the limit 
\[
	\alpha 
	:= \lim_{\Omega\ni y\to x_0}\uHp d_{x_0}(y)
\] 
exists. 
If $\alpha=0$, then $x_0$ must be regular by Theorem~\ref{thm:reg}, 
which is a contradiction. 
Hence $\alpha>0$.

\smallskip

\ref{d-lim} $\imp$ \ref{d-liminf} $\imp$ 
\ref{not-reg} $\imp$ \ref{not-reg-one-seq}
These implications are trivial.

\smallskip

$\neg$\ref{R} $\imp$ $\neg$\ref{not-reg-one-seq} 
Suppose that $x_0\in\overline{\{x\in\bdy\Omega:x\textup{ is regular}\}}$. 
For each integer $j\geq 2$, 
there exists a regular point $x_j\in B(x_0,1/j)\cap\bdy\Omega$. 
Define $f_j\in C(\bdy\Omega)$ 
by letting 
\[
	f_j(x) 
	= (jd_{x_0}(x)-1)_\limplus,
	\quad j=2,3,\ldots.
\]
Because $x_j$ is regular, 
there is $y_j\in B(x_j,1/j)\cap\Omega$ such that 
\[
	|\uHp f_j(y_j)| 
	= |f_j(x_j)-\uHp f_j(y_j)| 
	< 1/j.
\] 
Hence $y_j\to x_0$ and $\uHp f_j(y_j)\to 0$ as $j\to\infty$. 

Let $f\in C(\bdy\Omega)$ and $\alpha:=f(x_0)$. 
Let $\eps>0$. 
Then we can find an integer $k\geq 2$ such that 
$|f-\alpha|\leq\eps$ on $B(x_0,2/k)\cap\bdy\Omega$. 
Choose $m$ such that $|f-\alpha|\leq m$.
It follows that $f-\alpha\leq mf_j+\eps$ for every $j\geq k$, and thus 
\[
	\limsup_{j\to\infty}\uHp f(y_j)
	\leq \limsup_{j\to\infty}\uHp(mf_j+\alpha+\eps)(y_j)
	= m\lim_{j\to\infty}\uHp f_j(y_j)+\alpha+\eps
	= \alpha+\eps.
\]
Letting $\eps\to 0$ shows that 
$\limsup_{j\to\infty}\uHp f(y_j)\leq\alpha$. 

Applying this to $\ft=-f$ yields 
$\limsup_{j\to\infty}\uHp \ft(y_j)\leq-\alpha$. 
It follows that 
\[
	\liminf_{j\to\infty}\uHp f(y_j)
	\geq \liminf_{j\to\infty}\lHp f(y_j)
	= -\limsup_{j\to\infty}\uHp \ft(y_j)
	\geq \alpha,
\]
and hence 
$\lim_{j\to\infty}\uHp f(y_j)=f(x_0)$. 

\smallskip

\ref{Cp-V-bdy} $\eqv$ \ref{semireg-local}
Observe that \ref{Cp-V-bdy} is equivalent to 
the existence of a neighbourhood $U$ of $x_0$ 
with $\Cp(U\cap\bdy G)=0$, 
which is equivalent to \ref{semireg-local}, 
by the already proved equivalence \ref{Cp-V-bdy} $\eqv$ \ref{semireg} 
applied to $G$ instead of $\Omega$. 

\smallskip

\ref{Cp-V-bdy} $\imp$ \ref{Cp-V}
Let $V$ be a neighbourhood of $x_0$ such that $\Cp(V\cap\bdy\Omega)=0$. 
By Theorem~\ref{thm:irr-char-V}, 
\ref{V-Cp-V-bdy} $\imp$ \ref{V-alt-def-irr-super}, 
the set $U:=\Omega\cup(V\cap\bdy\Omega)$ is open and 
$\Cp(U\setm\Omega)=0$.

\smallskip

\ref{Cp-V} $\imp$ \ref{Cp-V-bdy}
This is trivial.

\smallskip

\ref{Cp-V} $\eqv$ \ref{alt-def-irr} $\eqv$ 
\ref{alt-def-irr-super}
In all three statements 
it follows directly that $V\subset\clOmega$. 
Thus their equivalence follows directly from Theorem~\ref{thm:irr-char-V}, 
with $V$ in Theorem~\ref{thm:irr-char-V} corresponding to 
$V\cap\bdy\Omega$ here.

\smallskip

\ref{alt-def-irr} $\imp$ \ref{rem-irr}
We only have to show the last part, i.e., 
that $x_0$ is irregular, 
but this follows from the already proved implication 
\ref{alt-def-irr} $\imp$ \ref{semireg}.

\smallskip

\ref{rem-irr} $\imp$ \ref{semireg} 
Let $f\in C(\bdy\Omega)$. 
Then $\uHp f$ has a \p-harmonic extension $U$ to $\Omega\cup V$ 
for some neighbourhood $V$ of $x_0$, 
and hence 
\[
	\lim_{\Omega\ni y\to x_0}\uHp f(y)
	= U(x_0).
\]
Since $x_0$ is irregular it follows that $x_0$ must be semiregular.

\smallskip

\ref{alt-def-irr-super} $\imp$ \ref{weaksemibarrier}
Let $u$ be a positive superharmonic function on $\Omega$. 
Then $\min\{u,1\}$ is superharmonic 
by Lemma~9.3 in Bj\"orn--Bj\"orn~\cite{BBbook}, 
and hence has a superharmonic extension $U$ to $\Omega\cup V$. 
As $U$ is lsc-regularized 
(see Section~\ref{sec:p-harmonic}) and 
$\mu(V\setm\Omega)=0$, 
it follows that $U\geq 0$ in $\Omega\cup V$. 
Suppose that $U(x_0)=0$. 
Then the strong minimum principle 
\cite[Theorem~9.13]{BBbook} 
implies that $U\equiv 0$ 
in the component of $\Omega\cup V$ that contains $x_0$. 
But this is in contradiction with $u$ being positive in $\Omega$, 
and thus 
\[
	\liminf_{\Omega\ni y\to x_0}u(x_0)
	\geq U(x_0)
	> 0.
\]
Thus there is no weak semibarrier at $x_0$.

\smallskip

$\neg$\ref{semibarrier} $\imp$ $\neg$\ref{weaksemibarrier} 
Let $u$ be a semibarrier at $x_0$. 
If $u>0$ in all of $\Omega$, 
then $u$ is a weak semibarrier at $x_0$. 
On the other hand, 
assume that there exists $x\in\Omega$ 
such that $u(x)=0$ 
(in this case $u$ is not a weak semibarrier). 
Then the strong minimum principle 
\cite[Theorem~9.13]{BBbook} 
implies that 
$u\equiv 0$ in the component $G\subset\Omega$ 
that contains $x$, 
and hence $x_0$ must be the only boundary point of $G$, 
because $u$ is a semibarrier. 
As $\Cp(X\setm G)\geq\Cp(X\setm\Omega)>0$, 
Lemma~4.3 in Bj\"orn--Bj\"orn~\cite{BB} 
(or Lemma~4.5 in \cite{BBbook}) 
implies that 
$\Cp(\{x_0\})=\Cp(\bdy G)>0$. 
By the Kellogg property (Theorem~\ref{thm:kellogg}), 
$x_0$ is regular, 
and hence 
Theorem~\ref{thm:reg} 
asserts that there is a positive 
barrier $v$ at $x_0$, 
and thus $v$ is a weak semibarrier. 

\smallskip

\ref{semibarrier} $\imp$ 
\ref{obst-dist-semibarrier}
This is trivial.

\smallskip

$\neg$\ref{R} $\imp$ $\neg$\ref{obst-dist-semibarrier}
Let $u$ be the continuous solution 
of the $\K_{d_{x_0},d_{x_0}}$-obstacle problem, 
which is superharmonic 
(see Section~\ref{sec:p-harmonic}). 
Moreover, it is clear that 
\[
	\liminf_{\Omega\ni y\to x}u(y)
	> 0
	\quad\text{whenever }x\in\bdy\Omega\setm\{x_0\},
\] 
and thus $u$ satisfies \ref{semibarrier-i} and \ref{semibarrier-iii} 
in Definition~\ref{def:semibarrier}.

Let $\{x_j\}_{j=1}^\infty$ be a sequence of regular boundary points 
such that $d_{x_0}(x_j)<1/j$. 
By Theorem~\ref{thm:reg}, 
$\lim_{\Omega\ni y\to x_j}u(y)=d_{x_0}(x_j)$. 
Hence we can find $y_j\in B(x_j,1/j)\cap\Omega$ so that 
$u(y_j)<2/j$. 
Thus $u$ satisfies 
\ref{semibarrier-ii} in Definition~\ref{def:semibarrier} as 
\[
	0 
	\leq \liminf_{\Omega\ni y\to x_0}u(y)
	\leq \liminf_{j\to\infty}u(y_j)
	= 0.
\]

\bigskip

We now assume that $\Omega$ is bounded or \p-parabolic. 

\smallskip

\ref{R} $\eqv$ \ref{rem-motiv} 
This follows directly from Theorem~\ref{thm:irr-char-V}, 
with $V$ therein 
corresponding to $V\cap\bdy\Omega$ here.
\end{proof}
We conclude our description of boundary points with some 
characterizations 
of strongly irregular points. 
As for regular and semiregular points, 
strong irregularity is a local property.
\begin{theorem}\label{thm:ess-irr-char}
Let $x_0\in\bdyOmegaX$\textup{,} $\delta>0$\textup{,} and 
$d_{x_0}\colon X^*\to[0,1]$ be defined by \eqref{eq-dx0}. 
Then the following are equivalent\/\textup{:} 
\begin{enumerate}
\item\label{ess-irr}
The point $x_0$ is strongly irregular.
\item\label{ess-local}
The point $x_0$ is strongly irregular with respect to 
$G:=\Omega\cap B(x_0,\delta)$.
\item\label{ess-one-seq}
The point $x_0$ is irregular and there exists a sequence 
$\{y_j\}_{j=1}^\infty$ in $\Omega$ such that 
$y_j\to x_0$ as $j\to\infty$\textup{,} 
and 
\[
	\lim_{j\to\infty}\uHp f(y_j)
	= f(x_0)
	\quad\text{for all }f\in C(\bdy\Omega).
\]
\item\label{ess-R}
It is true that $x_0\in\itoverline{R}\setm R$\textup{,} 
where $R:=\{x\in\bdy\Omega:x\text{ is regular}\}$. 
\item\label{ess-d-liminf} 
It is true that 
\[
	\liminf_{\Omega\ni y\to x_0}\uHp d_{x_0}(y)
	= 0
	< \limsup_{\Omega\ni y\to x_0}\uHp d_{x_0}(y).
\]
\item\label{ess-f-nolim} 
There exists $f\in C(\bdy\Omega)$ such that 
\[
	\lim_{\Omega\ni y\to x_0}\uHp f(y)
\]
does not exist. 
\item\label{obst-dist-irr} 
The continuous solution $u$ of the $\K_{d_{x_0},d_{x_0}}$-obstacle problem  satisfies 
\[
	\liminf_{\Omega\ni y\to x_0}u(y)
	= 0
	< \limsup_{\Omega\ni y\to x_0}u(y).
\]
\item\label{barrier-irr}
There is a semibarrier \textup{(}or equivalently there is 
a weak semibarrier\/\textup{)} but no barrier at $x_0$. 
\end{enumerate}
\end{theorem}
The trichotomy property (Theorem~\ref{thm:trichotomy}) shows that 
a boundary point 
is either regular, semiregular, or strongly 
irregular. 
We will use this in the following proof.
\begin{proof}
\ref{ess-irr} $\eqv$ \ref{ess-local}
By Theorems~\ref{thm:reg} and~\ref{thm:rem-irr-char}, 
regularity and semiregularity are local properties, 
and hence this must be true also for strong irregularity. 

\smallskip

\ref{ess-irr} $\eqv$ \ref{ess-one-seq} $\eqv$ \ref{ess-R}
This follows from 
Theorem~\ref{thm:rem-irr-char} 
\ref{semireg} $\eqv$ \ref{not-reg-one-seq} $\eqv$ \ref{R}.

\smallskip

\ref{ess-irr} $\imp$ \ref{ess-d-liminf}
Since $x_0$ is strongly irregular and $\uHp d_{x_0}$ is nonnegative, 
it follows that 
\[
	\liminf_{\Omega\ni y\to x_0}\uHp d_{x_0}(y)
	= 0
	\leq \limsup_{\Omega\ni y\to x_0}\uHp d_{x_0}(y).
\]
If $\limsup_{\Omega\ni y\to x_0}\uHp d_{x_0}(y)=0$, 
then $x_0$ must be regular by Theorem~\ref{thm:reg}, 
which is a contradiction. 
Thus 
\[
	\limsup_{\Omega\ni y\to x_0}\uHp d_{x_0}(y)
	> 0.
\]

\smallskip

\ref{ess-d-liminf} $\imp$ \ref{ess-f-nolim}
This is trivial.

\smallskip

\ref{ess-f-nolim} $\imp$ \ref{ess-irr}
By definition, $x_0$ is neither regular nor semiregular, 
and hence must be strongly irregular.

\smallskip

\ref{ess-irr} $\eqv$ \ref{obst-dist-irr}
Theorem~\ref{thm:reg} 
shows that $x_0$ is regular if and only if 
$\lim_{\Omega\ni y\to x_0}u(y)=0$. 
On the other hand, 
Theorem~\ref{thm:rem-irr-char} implies that $x_0$ is semiregular 
if and only if 
$\liminf_{\Omega\ni y\to x_0}u(y)>0$. 
The equivalence follows by combining these two facts. 

\smallskip

\ref{ess-irr} $\eqv$ \ref{barrier-irr}
By Theorem~\ref{thm:rem-irr-char}, 
$x_0$ is semiregular if and only if there is no (weak) 
semibarrier at $x_0$. 
On the other hand, 
by Theorem~\ref{thm:reg}, 
there is a barrier at $x_0$ if and only if $x_0$ is regular. 
Combining these two facts gives the equivalence.
\end{proof}


\end{document}